\documentclass{article} 
\usepackage{theorem, graphicx}

{\theorembodyfont{\normalfont}

\newtheorem{example}{Example}

}

\begin{document}

\title{Logic of Simultaneity}            

\author{Kenji Tokuo\footnote{Department of Computer and Control Engineering, Oita National College of Technology, 1666 Maki, Oita 870-0152, Japan.  Email: tokuo@oita-ct.ac.jp}}
\date{}

\maketitle

\begin{abstract} 
A logical model of spatiotemporal structures is pictured as a succession of processes in time. 
One usual way to formalize time structure is to assume the global  existence of time points and then collect some of them to form time intervals of processes. 
Under this set-theoretic approach, the logic that governs the processes acquires a Boolean structure. 
However, in a real distributed system or a relativistic universe where the message-passing time between different locations is not negligible, the logic has no choice but to accept time interval instead of time point as a primitive concept. 
From this modeling process of spatiotemporal structures, orthologic, the most simplified version of quantum logic, emerges naturally. 
\end{abstract}

%\Keywords{simultaneity, relativity theory, temporal logic, orthologic, quantum logic}

\section{General Assumptions}
We consider a simple model of concurrent processes in computer systems, or more generally, a logical  representation of the physical universe that consists of causal sequences of events occurring at different spatial locations. These locations are here referred to as {\em sites}\/. Formally, a site is a sequence of processes ordered by the happened-before relation. 
The term {\em process}\/ denotes one of the continuing states in a site, including in particular nothing happening states. It is also assumed for simplicity that (i) the processes in each site occur consecutively with no time gaps, (ii) each process lasts for some non-zero time duration. Any change of processes is referred to as an {\em event}\/ (Fig. 1).  

\begin{center}
\includegraphics[scale=0.6]{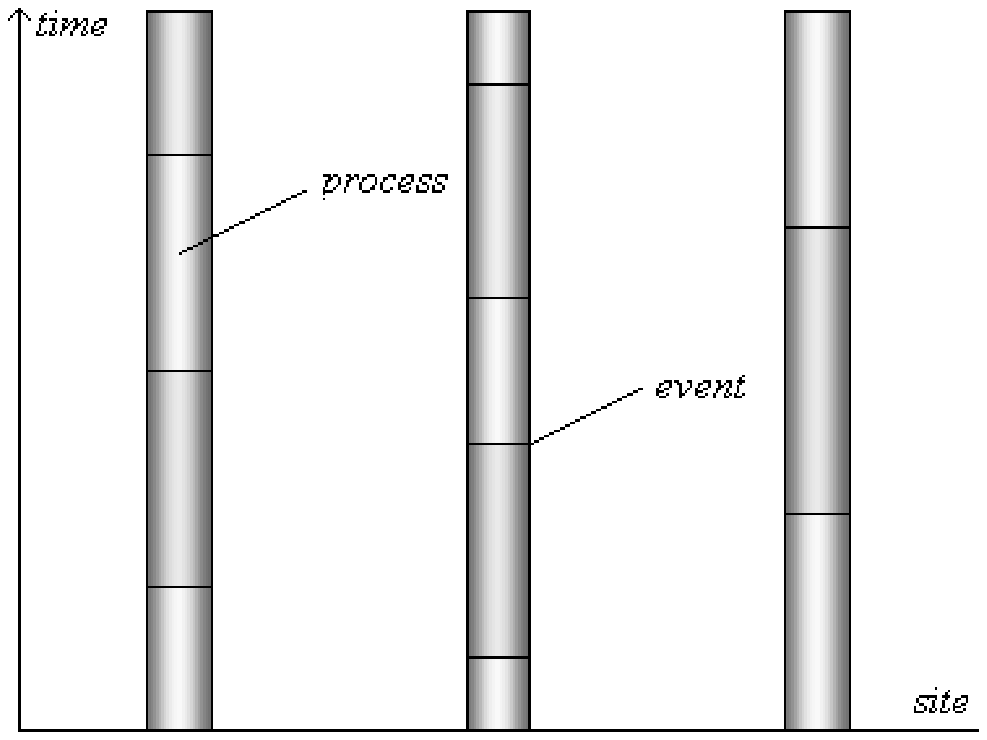}\\

[Fig. 1]
\end{center}

\section{Concept of Time in Non-Relativistic Theories}
For the later comparison, this section is devoted to introduce a formalization of time in non-relativistic situations. The earliest such treatment is due to Russel (1926).

\subsection{Logical Construction of Time}
Since we assume that any process is not instantaneous, but occupies some non-zero time duration, we can say that a process \(p\) {\em is earlier than}\/ a process \(q\) if \(p\) ends before \(q\) begins, and that \(p\) {\em is simultaneous with}\/ \(q\) if \(p\) partly or completely overlaps with \(q\), i.e. neither \(p\) is earlier than \(q\) nor vice versa. In particular, \(p\) is simultaneous with \(p\) itself. Note that this relation of simultaneity is reflexive and symmetric, but not transitive. 

A maximal set of simultaneity, which is referred to as a {\em time point}, is then defined as a maximal set of processes (with respect to the set inclusion ordering), any two members of which are simultaneous with each other. The collection of all time points is denoted by \(\cal P\), i.e.

\begin{equation}
{\cal P} \equiv \{ T \in 2^{\bf Proc} \,\,|\,\, \forall p\in T.\, S(p, q) \Leftrightarrow q \in T \} 
\end{equation}

\noindent
where {\bf Proc} denotes the set of all processes, \(2^{\bf Proc}\) denotes the power set of {\bf Proc} (the set of all subsets of {\bf Proc}), and \(S\) denotes the simultaneity relation. 

\begin{example}
In Fig. 2, we have
\(T_{1} = \{ p_1, q_1, r_1 \}\), \(T_2 = \{ p_1, q_2, r_1 \}\), \(T_3 = \)
\(\{ p_2, q_2, r_1 \}\), \(T_4 = \{ p_2, q_2, r_2 \}\), \(T_5 = \{ p_2, q_3, r_2 \}\), \(T_6 = \{ p_3, q_3, r_2 \}\), \(T_7 = \{ p_3, q_4,\) \(r_2 \}\), \(T_8 = \{ p_3, q_4, r_3 \}\), \(T_9 = \{ p_4, q_4, r_3 \}\), \(T_{10} = \{ p_4, q_5, r_3 \}\).
The broken lines represent the simultaneous time points.
\begin{center}
\includegraphics[scale=0.6]{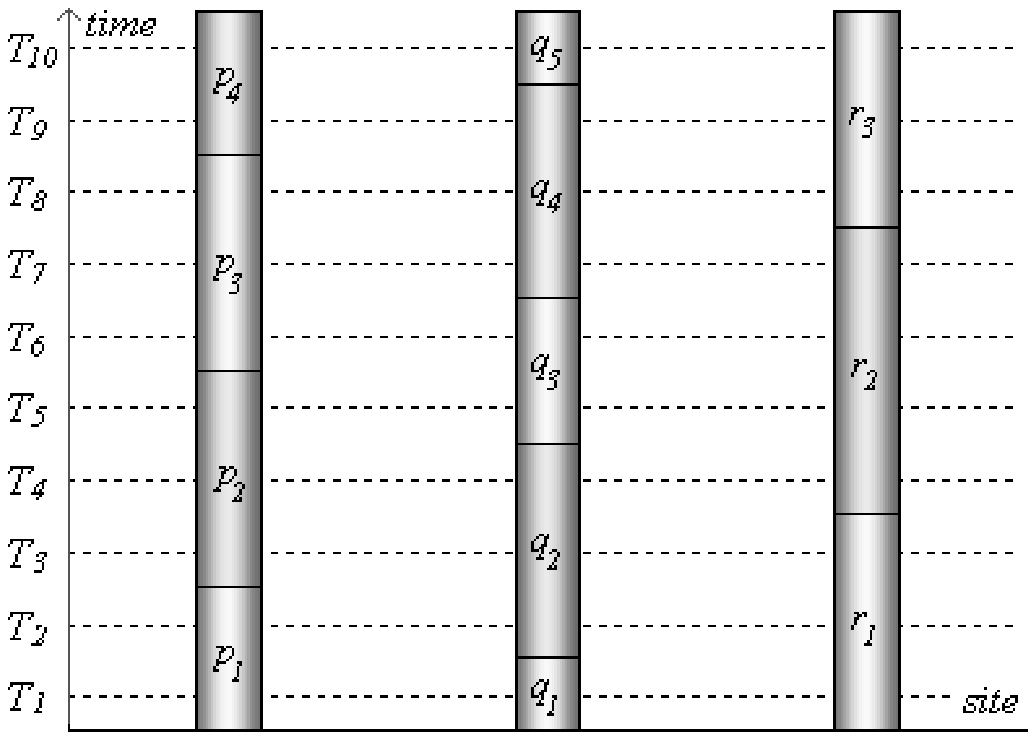}\\

[Fig. 2]
\end{center}
\end{example}

It is worth noting that the happened-before relation on {\bf Proc} defines a linear ordering on \(\cal P\). For any distinct time points \(T\) and \(T'\), there must be a process \(p\) that is in \(T\) but not in \(T'\), and a process \(p'\) that is in \(T'\) but not in \(T\). By the definition of time point, either of the following holds: ``\(p\) is earlier than \(p'\)'' or ``\(p'\) is earlier than \(p\).''  In the former case let \(T < T'\) and in the latter case \(T' < T\). 

Time points are used to introduce the concept of {\em time intervals.}\/ For each process \(p\), its time interval \([\, p\, ]\) is defined as the set of all time points that contain \(p\), i.e. 

\begin{equation}
[\, p\, ] \equiv \{ T \in {\cal P} \,\,|\,\, p \in T \}\textrm{.}
\end{equation}

\begin{example}
In Fig. 2, we see that \(T_1 < T_2 < \cdots. < T_{10}\), and that \([\, p_1\, ] = \{ T_1, T_2 \}\), \([\, q_1\, ] = \{ T_1\}\),  \([\, r_1\, ] = \{ T_1, T_2, T_3 \}, \ldots \).   
\end{example}

\subsection{Logic of Simultaneity in Non-Relativistic Theories}
We are now led to the logic of simultaneity by considering a process \(p\) and its time interval \([\, p\, ]\) as an atomic {\em proposition}\/ and its {\em truth value}, respectively. An atomic proposition \(p\) asserts that {\em the process \(p\) is occurring.}\/ Complex propositions are built out of atomic propositions and logical connectives: 

\begin{itemize}
\item For any proposition \(p\), \(\neg p\) denotes the proposition that the proposition \(p\) is not true. The truth value of \(\neg p\) is defined as \([\, \neg p\, ] \equiv {\cal P} - [\, p\, ]\) (set-theoretic complement relative to \(\cal P\)), which amounts to the time interval that \(p\) is not true. 
\item For any propositions \(p\) and \(q\), \(p \wedge q\) denotes the proposition that the propositions \(p\) and \(q\) are both true. The truth value of \(p \wedge q\) is defined as \([\, p \wedge q\, ] \equiv [\, p\, ] \cap [\, q\, ]\) (set-theoretic intersection), which amounts to the time interval that \(p\) and \(q\) are both true. 
\item For any propositions \(p\) and \(q\), \(p \vee q\) denotes the proposition that at least one of the propositions \(p\) and \(q\) is true. The truth value of \(p \vee q\) is defined as \([\, p \vee q\, ] \equiv [\, p\, ] \cup [\, q\, ]\) (set-theoretic union), which amounts to the time interval that at least one of \(p\) and \(q\) is true. 
\end{itemize}

\noindent
Since the operations coincide with the usual set-theoretic ones, it is obvious that the resulting logic is Boolean. 

\section{Concept of Time in Relativistic Theories}

\subsection{Undecidability  of  Simultaneity in Relativistic Theories}

In the preceding section, we have implicitly assumed the existence of the global clock, which is represented by the linearly arranged time points. However, the classical concept of simultaneity loses its meaning in a real distributed system or a relativistic universe. What the principle of special relativity says is that it does take a non-zero time duration to transmit any causal signals between spatially separated sites (For a basic reference, see Taylor (1992)). As shown in Fig. 3, we refer to any signal capable of transmitting information between sites as a {\em message}. 

\begin{center}
\includegraphics[scale=0.6]{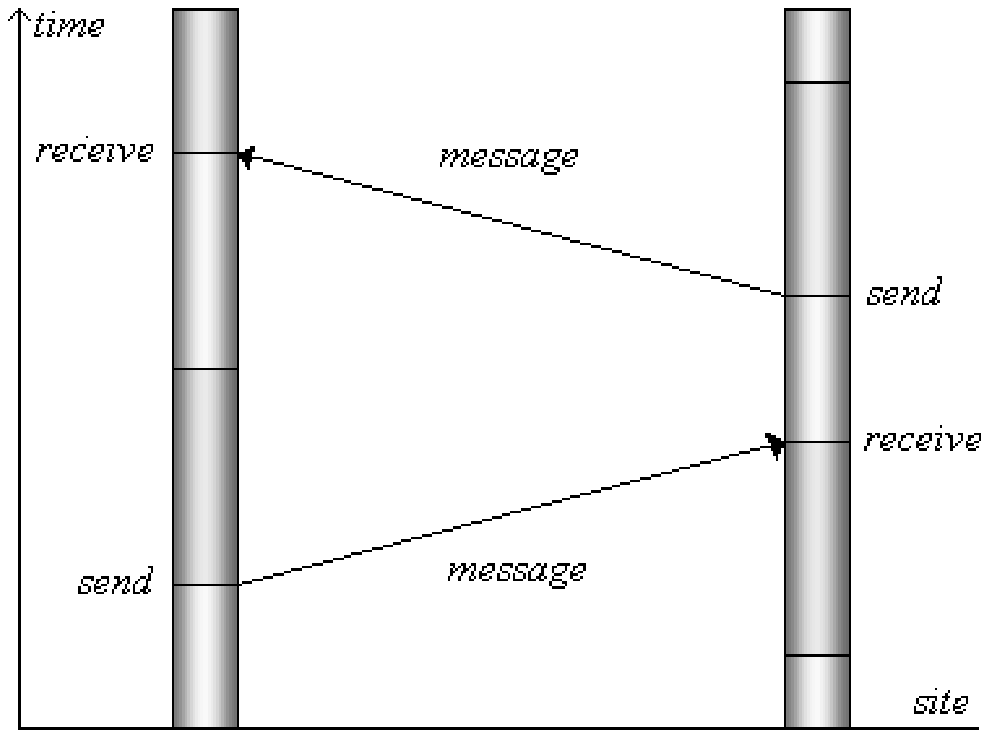}\\

[Fig. 3]
\end{center}

While a message transmission is needed to synchronize two clocks at different sites, 
the transmission time duration is not measurable without reference to synchronized clocks located at both sites; it is a vicious circle. Hence we must abandon the attempt to provide the global time points, i.e. the clock common to all sites. 
Thus the simultaneity based on the overlap relation cannot be defined since it cannot generally be determined whether two processes at different sites overlap temporally.

\begin{example}
In a relativistic universe, it is essentially meaningless to distinguish between A and B (Fig. 4). 
\begin{center}
\includegraphics[scale=0.6]{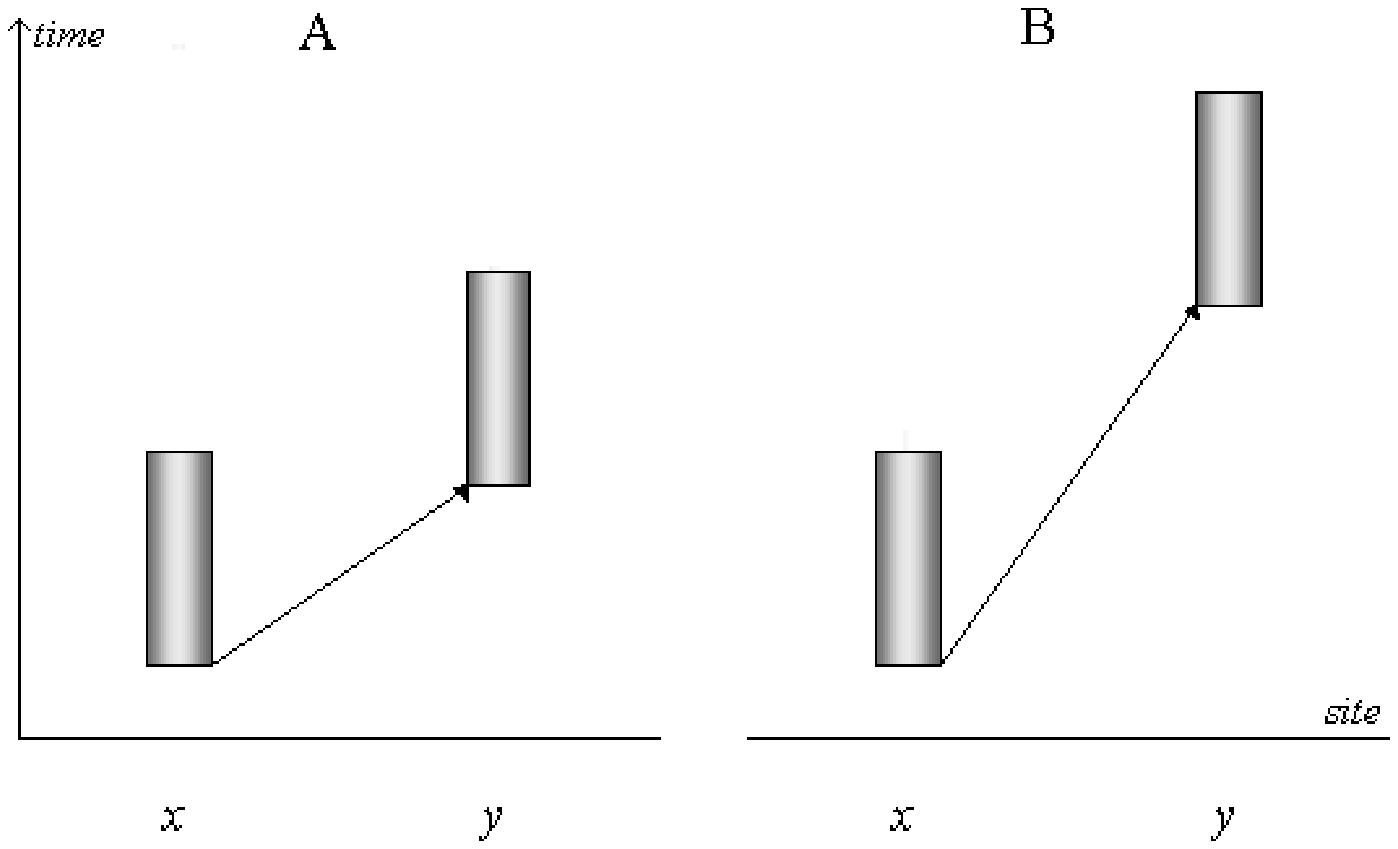}\\

[Fig. 4]
\end{center}
\end{example}

\subsection{Temporal Containment}
We should therefore focus on the special cases where we can say with certainty that two or more processes run simultaneously. In simple cases where a message is sent from a site \(x\) to a site \(y\) and then another message is sent back from \(y\) to \(x\), it is verifiable that the processes in \(x\) that occur consecutively with no time gaps between the sending event and the receiving event {\em temporally contain}\/ the processes in \(y\) that occur between the receiving event and the sending event. 
\begin{example}
Temporal containment: a process \(q\) is occurring whenever a process \(r\) is occurring (Fig. 5).
\begin{center}
\includegraphics[scale=0.6]{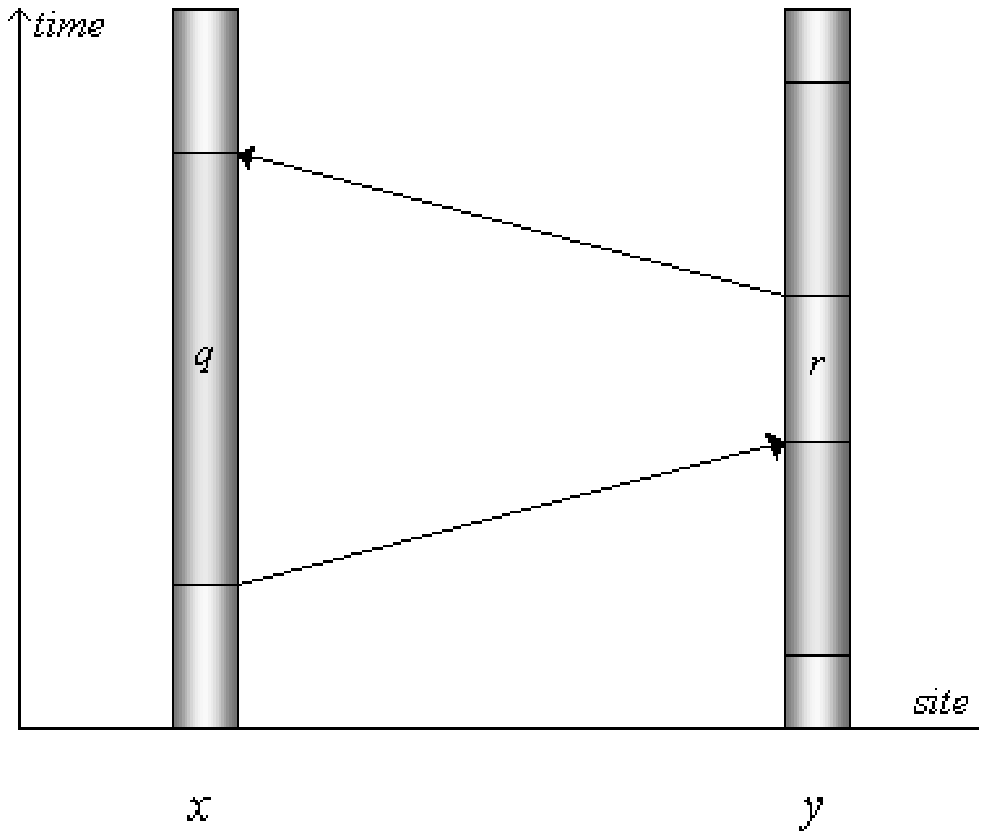}\\

[Fig. 5]
\end{center}
\end{example}

To sum up, a relativistic universe admits simultaneity only in the sense of temporal containment shown in Fig 5. In the following discussion, we drop the idea of defining time interval via simultaneous time points, but directly construct the logic employing temporal containment relation between processes. 

\subsection{Logic of Simultaneity in Relativistic Theories}
We say that two processes have a {\em causal relationship}\/ in a relativistic universe if they are linked with the happened-before relation (Lamport, 1978). The happened-before relation \(B\) on {\bf Proc} is defined as the smallest relation satisfying the following conditions: (i) If \(p\) and \(q\) are processes of the same site, and \(p\) occurs before \(q\), then \(B(p, q)\). (ii) If \(p\) ends with an event of sending a message and \(q\) begins with an event of receiving that message, then \(B(p, q)\). (iii) If \(B(p, q)\) and \(B(q, r)\), then \(B(p, r)\). The causality relation \(C\) on {\bf Proc} is then defined as the smallest relation satisfying the following condition: \(C(p, q)\) if and only if \(B(p, q)\) or \(B(q, p)\). Note that the relation \(C\) is irreflexive and symmetric, but not transitive. 

Using the same notation as before, we informally denote by \([\, p\, ]\) the time interval of a process \(p\). The example shown in Fig. 5 thus indicates that \([\, q\, ] \supseteq [\, r\, ]\). This containment relation is characterized by the fact that {\em any process that has a causal relationship with \(q\) has a causal relationship with \(r\)}, i.e. 

\begin{equation}\label{eq:PCC}
\forall p.\, (C(p, q) \Rightarrow C(p, r) )
\end{equation}

To formalize a more general setting where a process is covered by two or more processes (Fig. 6), we need a slight modification of the formula (\ref{eq:PCC}).  Letting \([\, q_1, q_2, \ldots\, ]\) be the time interval that at least one of \(\{ q_1, q_2, \ldots \}\) occurs, we say that \([\, q_1, q_2, \ldots\, ] \supseteq [\, r\, ]\) if {\em any process that has a causal relationship with any process of \(\{ q_1, q_2, \ldots\, \}\) has a causal relationship with \(r\)}, i.e.

\begin{equation}\label{eq:CC}
\forall p.\, ( (\forall q \in \{ q_1, q_2, \ldots \}.\, (C(p, q))) \Rightarrow C(p, r)) 
\end{equation}

\begin{center}
\includegraphics[scale=0.6]{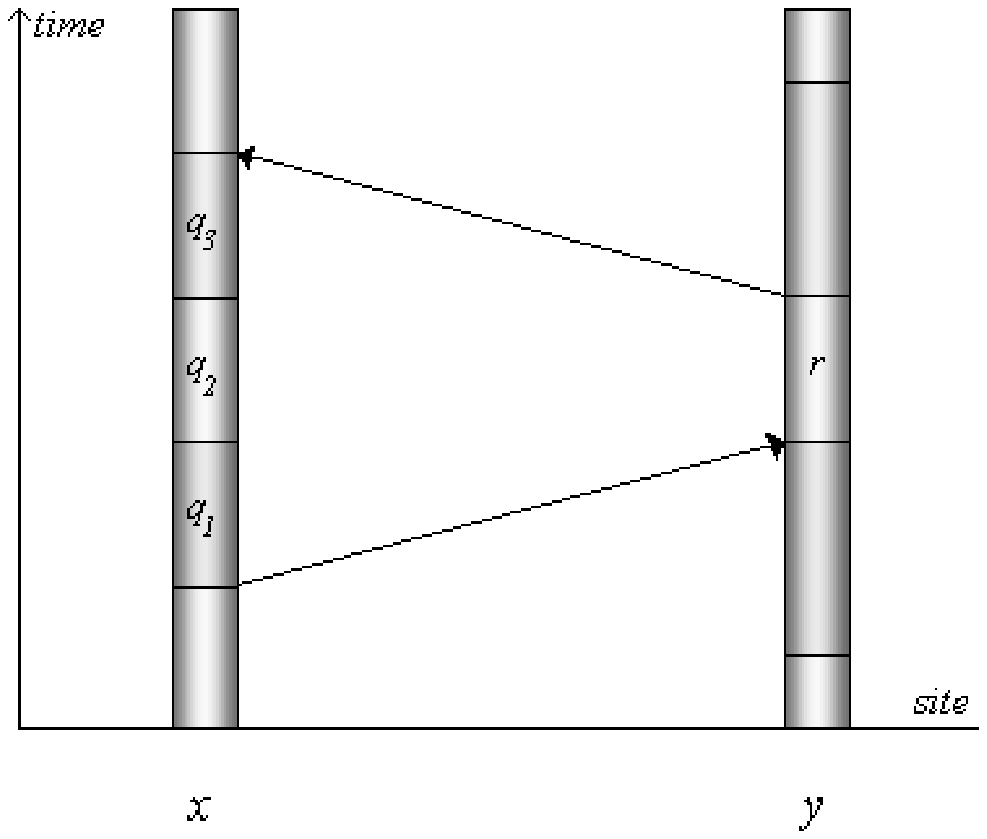}\\

[Fig. 6]
\end{center}

Now taking the containment relation as fundamental, we can conceive of a non-Boolean model for the logic of  spatiotemporal structures. 
Letting 
\[{\cal I} \equiv \{ I\in 2^{\bf Proc} \,\,|\,\, \forall r.\, ( \forall p.\, (\forall q\in I.\, C(p, q) \Rightarrow C(p, r) ) \Leftrightarrow r \in I) \}\]
be the collection of all time intervals, we stipulate that the time interval \([\, q_1, q_2, \ldots\, ]\) be the minimum element (with respect to the set inclusion ordering) in \(\cal I\) that contains the set \(\{ q_1, q_2, \ldots \}\). This definition of time interval indeed satisfies the condition (\ref{eq:CC}). The existence of the minimum element is assured by the fact that \(\cal I\) is closed under the set-theoretic intersection. Note that \(\cal I\) is {\em not}\/ closed under the set-theoretic union. 

\begin{example}
In Fig. 7, we have
\small{\({\cal I} = \{ \phi\), \(\{p_1\}\), \(\{p_2\}\), \(\{p_1, p_2\}\), \(\{p_3\}\), \(\{p_1, p_3\}\), 
\(\{p_4\}\), \(\{p_1, p_4\}\), \(\{p_2, p_4\}\), \(\{p_1, p_2, p_4\}\), \(\{p_3, p_4\}\), \(\{p_1, p_3, p_4\}\), \(\{q_1\}\), \(\{q_2\}\), \(\{p_1, q_1, \)
\(q_2\}\), \(\{q_3\}\), \(\{p_2, p_3, q_3\}\), \(\{q_1, q_3\}\), \(\{q_2, q_3\}\), \(\{p_1, q_1, q_2, q_3\}\), \(\{p_1, p_2, p_3, q_1, q_2, q_3\}\), 
\(\{q_4\}\), \(\{q_1, q_4\}\), \(\{q_2, q_4\}\), \(\{p_1, q_1, q_2, q_4\}\), \(\{q_3, q_4\}\), \(\{q_1, q_3, q_4\}\), \(\{q_5\}\), \(\{q_1, q_5\}\), \(\{q_2, \)
\(q_5\}\), \(\{p_1, q_1, q_2, q_5\}\), \(\{q_3, q_5\}\), \(\{q_1, q_3, q_5\}\), \(\{q_2, q_3, q_5\}\), \(\{p_1, q_1, q_2, q_3, q_5\}\), \(\{p_4, q_4, \)
\(q_5\}\), \(\{p_4, q_1, q_4, q_5\}\), \(\{p_4, q_2, q_4, q_5\}\), \(\{p_1, p_4, q_1, q_2, q_4, q_5\}\), \(\{p_4, q_3, q_4, q_5\}\), \(\{p_2, p_3,  \)
\(p_4, q_3, q_4, q_5\}\), \(\{p_4,\) \(q_1, q_3, q_4, q_5\}\), \(\{q_1, r_1\}\), \(\{r_2\}\), \(\{q_2, q_3, q_4, r_2\}\), \(\{q_1, r_1, r_2\}\), \(\{p_1, \)
\(q_1, q_2, q_3, q_4, r_1, r_2\}\), \(\{q_5, r_3\}\), \(\{q_1, q_5, r_1,\) \(r_3\}\), \(\{q_5, r_2, r_3\}\), \({\bf Proc} \}\)}.
\begin{center}
\includegraphics[scale=0.6]{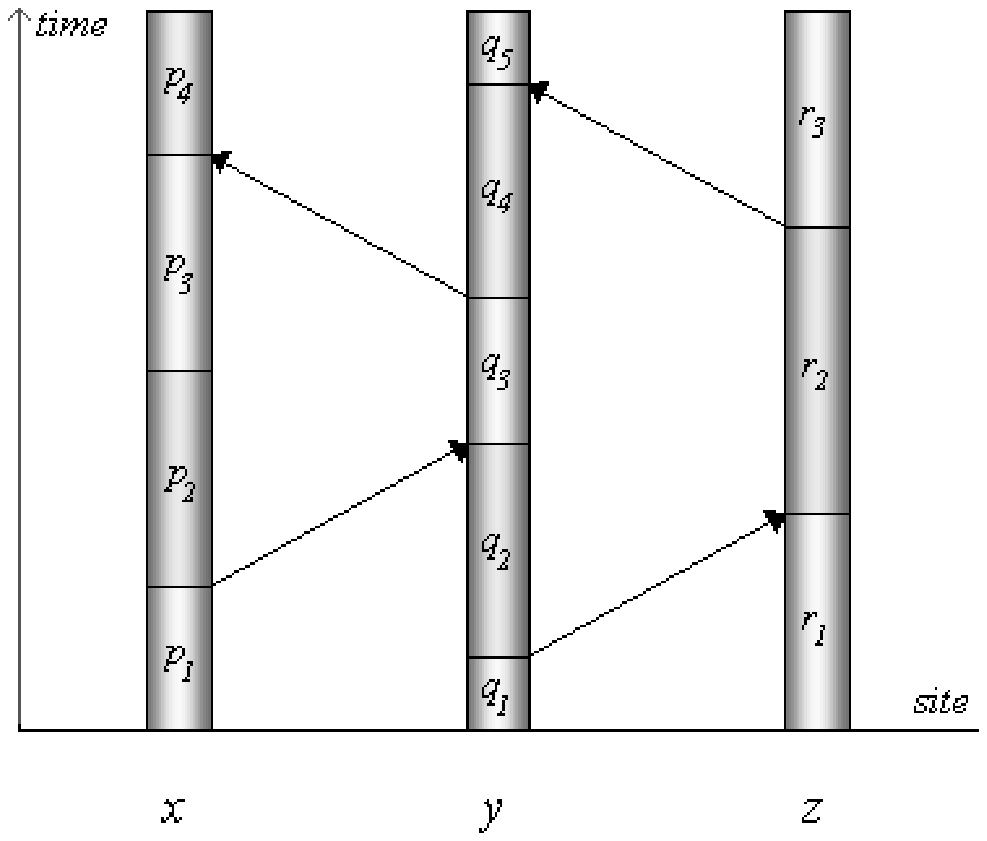}\\

[Fig. 7]
\end{center}
\end{example}

As in the non-relativistic case, we obtain a logic from time structure by considering the time interval of a process as the truth value of an atomic proposition. Now the resulting algebra is {\em ortholattice}, which is not Boolean in general. The associated logic is called {\em orthologic} or  {\em minimal quantum logic} (For a basic reference, see Birkhoff (1967) and Dalla Chiara (2001),): 

\begin{itemize}
\item For any proposition \(p\), \(\neg p\) denotes the proposition that the proposition \(p\) is not true. The truth value of \(\neg p\) is defined as
\([\, \neg p\, ] \equiv \{ q \,\,|\,\, \forall r\in [\, p\, ].\, C(q, r)\}\)
(orthocomplement relative to \(\cal I\)), which \\
amounts to the time interval that \(p\) is not true. The following facts follow from the definition. 
\begin{itemize}
\item \([\, \neg p\, ] \in {\cal I}\)
\item \([\, \neg\neg p\, ] = [\, p\, ]\) 
\item \([\, p\, ] \subseteq [\, q\, ]\) if and only if \([\, \neg q\, ] \subseteq [\, \neg p\, ]\) 
\end{itemize}
\item For any propositions \(p\) and \(q\), \(p \wedge q\) denotes the proposition that the propositions \(p\) and \(q\) are both true. The truth value of \(p \wedge q\) is defined as \([\, p \wedge q\, ] \equiv [\, p\, ] \cap [\, q\, ]\) (set-theoretic intersection), which amounts to the time interval that \(p\) and \(q\) are both true. Since \(\cal I\) is closed under the set-theoretic intersection, \(\wedge\) corresponds to the infimum operator on \(\cal I\) with respect to the set inclusion ordering. 
\item For any propositions \(p\) and \(q\), \(p \vee q\) denotes the proposition that at least one of the propositions \(p\) and \(q\) is true. The truth value of \(p \vee q\) is defined as \([\, p \vee q\, ] \equiv [\, \neg(\neg p \wedge \neg q)\, ]\), which amounts to the time interval that at least one of \(p\) and \(q\) is true. Since \(\neg\) has the above-mentioned properties and \(\wedge\) is the infimum operator on \(\cal I\), \(\vee\) corresponds to the supremum operator on \(\cal I\) with respect to the set inclusion ordering. Note that we have
\([\, p_1, p_2, \ldots\, ] = [\, p_1 \vee p_2 \vee \ldots\, ]\) 
for atomic propositions \(p_1, p_2, \ldots\).
\end{itemize}

\begin{example}
A typical statement which is always true in Boolean logic but not necessarily true in orthologic is the distributive law of \(\wedge\) over \(\vee\), i.e.
\begin{equation}
[\, (p \vee q) \wedge r\, ] = [\, (p \wedge r) \vee (q \wedge r)\, ]. 
\end{equation}
In Fig. 7, we can find a counterexample to distributivity: since \([\, p_2 \vee p_3\, ] = \{ p_2, p_3, q_3 \}\) and \([\, q_3\, ] = \{ q_3 \}\), we infer \([\, (p_2 \vee p_3) \wedge q_3\, ] = \{ q_3 \}\), while since \([\, p_2 \wedge q_3\, ] = \phi\) and \([\, p_3 \wedge q_3\, ] = \phi\), we infer \([\, (p_2 \wedge q_3) \vee (p_3 \wedge q_3)\, ] = \phi\). The failure of distributivity illustrates the fact that the analysis of spatiotemporal structures deduces non-Boolean orthologic when global synchronized clocks are not available. 
\end{example}

\bibliographystyle{plain}

%\end{paper}
\end{document}